\documentclass[12pt]{amsart}
\usepackage{graphicx,amsmath,amssymb}

\newcommand{\C}{\mathbb C}
\newcommand{\R}{\mathbb R}
\newcommand{\Z}{\mathbb Z}
\newcommand{\N}{\mathbb N}
\renewcommand{\d}{\prime} 
\newcommand{\dd}{{\prime \prime}}

\renewcommand{\Re}{{\rm Re}\,}
\renewcommand{\Im}{{\rm Im}\,}
\newtheorem{theorem}{Theorem}
\newtheorem{lemma}[theorem]{Lemma}
 
\newtheorem{corollary}[theorem]{Corollary}

\newtheorem{definition}{Definition}
\newtheorem*{remark}{Remark}

\begin{document}
\title[]
{Eigenvalues of $\mathcal{PT}$-symmetric oscillators with polynomial potentials}
\author[]
{Kwang C. Shin}
\address{Department of Mathematics, University of Missouri, Columbia, MO 65211, USA}
\date{January 02, 2005}

\begin{abstract}
We study the eigenvalue problem $-u^\dd(z)-[(iz)^m+P_{m-1}(iz)]u(z)=\lambda u(z)$ with the boundary conditions that $u(z)$ decays to zero as $z$ tends to infinity along the rays $\arg z=-\frac{\pi}{2}\pm \frac{2\pi}{m+2}$, where $P_{m-1}(z)=a_1 z^{m-1}+a_2 z^{m-2}+\cdots+a_{m-1} z$ is a polynomial and integers $m\geq 3$. We provide an asymptotic expansion of the eigenvalues $\lambda_n$ as $n\to+\infty$, and prove that for each {\it real} polynomial $P_{m-1}$, the eigenvalues are all real and positive, with only finitely many exceptions. 
\end{abstract}

\maketitle

\begin{center}
{\it 2000 {\it Mathematics subject classification}: 34L40, 34L20}
\end{center}

\baselineskip = 18pt

\section{Introduction}
\label{introduction}
For integers $m\geq 3$ fixed, we are considering the ``non-standard'' non-self-adjoint eigenvalue problems
\begin{equation}\label{ptsym}
H u(z,\lambda):=\left[-\frac{d^2}{dz^2}-(iz)^m-P_{m-1}(iz)\right]u(z,\lambda)=\lambda u(z,\lambda),\quad\text{for some $\lambda\in\C$},
\end{equation}
with the boundary condition that 
\begin{equation}\label{bdcond}
\text{$u(z,\lambda)\rightarrow 0$ exponentially, as $z\rightarrow \infty$ along the two rays}\quad \arg (z)=-\frac{\pi}{2}\pm \frac{2\pi}{m+2},
\end{equation}
where $P_{m-1}$ is a polynomial of degree at most $m-1$ of the form 
\begin{equation}\label{poly_form}
P_{m-1}(z)=a_1z^{m-1}+a_2z^{m-2}+\cdots+a_{m-1}z,\quad a_j\in\C\,\,\text{\,for $1\leq j\leq m-1$}.
\end{equation}
We let $$a:=(a_1,a_2,\ldots, a_{m-1})\in\C^{m-1}$$ be the coefficient vector of $P_{m-1}(z)$. We are mainly interested in the case when $P_{m-1}$ is real, that is, when $a\in\R^{m-1}$. However, some interesting facts in this paper hold also for $a\in\C^{m-1}$. So except for Theorem \ref{main_theorem} below, we will use $a\in\C^{m-1}$.

If a nonconstant function $u$ satisfies \eqref{ptsym} with some $\lambda\in\C$ and the boundary condition \eqref{bdcond}, then we call $\lambda$ an {\it eigenvalue} of $H$ and $u$ an {\it eigenfunction of $H$ associated with the eigenvalue $\lambda$}.  Also, the {\it geometric multiplicity of an eigenvalue $\lambda$} is the number of linearly independent eigenfunctions  associated with the eigenvalue $\lambda$. 
The  operator $H$ in \eqref{ptsym} with potential $V(z)=-(iz)^m-P_{m-1}(iz)$ is called {\it $\mathcal{PT}$-symmetric} if $\overline{V(-\overline{z})}=V(z)$, $z\in\C$. Note that $V(z)=-(iz)^m-P_{m-1}(iz)$ is a $\mathcal{PT}$-symmetric potential if and only if $a\in\R^{m-1}$. 


Before we state our main theorems, we first introduce some known facts by Sibuya \cite{Sibuya} about the eigenvalues $\lambda$ of $H$.
\begin{theorem}\label{main2}
The eigenvalues $\lambda_n$ of $H$ have the following properties.
\begin{enumerate}
\item[(I)] The set of all eigenvalues is a  discrete set in $\C$.
\item[(II)] The geometric multiplicity of every eigenvalue is one.
\item[(III)] Infinitely many eigenvalues, accumulating at infinity, exist.
\item[(IV)] The eigenvalues have the following asymptotic expansion
\begin{equation}\label{Bender_exp}
\lambda_n =\left(\frac{\Gamma\left(\frac{3}{2}+\frac{1}{m}\right)\sqrt{\pi}\left(n-\frac{1}{2}\right)}{\sin \left(\frac{\pi}{m}\right)\Gamma\left(1+\frac{1}{m}\right) } \right)^{\frac{2m}{m+2}}[1+o(1)]\quad\text{as $n$ tends to infinity},\quad n \in \N,
\end{equation}
where the error term $o(1)$ could be complex-valued.
\end{enumerate} 
\end{theorem}

This paper is organized as follows. In Section \ref{prop_sect}, we will introduce work of Hille \cite{Hille} and Sibuya \cite{Sibuya}, regarding properties of solutions of \eqref{ptsym}. We then improve on the asymptotics of a certain function in \cite{Sibuya}. In Section \ref{entire_sect}, we introduce an entire function $C(a,\lambda)$ whose zeros are the eigenvalues of $H$, due to Sibuya \cite{Sibuya}. In Section \ref{asymp_sect}, we then provide asymptotics of $C(a,\lambda)$ as $\lambda\to\infty$ in the complex plane, improving the  asymptotics of $C(a,\lambda)$ in \cite{Sibuya}. In Section \ref{asymp_eigen}, we will improve the asymptotic expansion \eqref{Bender_exp} of the eigenvalues. In particular, we will prove the following. Throughout this paper, we use that $\lfloor x\rfloor$ is the largest integer that is less than or equal to $x\in\R$.
\begin{theorem}\label{main_thm1}
Let $a\in\C^{m-1}$ be fixed. Then there exist $e_{\ell}(a)\in\C, \,\, 1\leq \ell\leq \frac{m}{2}+1$ such that the eigenvalues $\lambda_n$ of $H$ have the asymptotic expansion
\begin{equation}\label{asym_eq1}
\lambda_n\underset{n\to+\infty}{=}\lambda_{0,n}+\sum_{\ell=1}^{\lfloor\frac{m}{2}+1\rfloor}e_{\ell}(a)
\lambda_{0,n}^{1-\frac{\ell}{m}}+o\left(\lambda_{0,n}^{\frac{1}{2}-\frac{1}{m}}\right),
\end{equation}
where
\begin{equation}\nonumber
\lambda_{0,n}=\left(\frac{\left(n+\frac{1}{2}\right)\pi}{K_{m}\sin\left(\frac{2\pi}{m}\right)}\right)^{\frac{2m}{m+2}}\quad\text{with}\quad K_{m}=\int_0^{\infty}\left(\sqrt{1+t^m}-\sqrt{t^m}\right)\, dt>0.
\end{equation}
\end{theorem}
One can compute $K_m$ directly (or see equation (2.22) in \cite{Dorey2} with the identity $\Gamma(s)\Gamma(1-s)=\pi\csc(\pi s)$) and obtains
$$K_m=\frac{\sqrt{\pi}\Gamma\left(1+\frac{1}{m}\right)}{2\cos\left(\frac{\pi}{m}\right)\Gamma\left(\frac{3}{2}+\frac{1}{m}\right)}.$$
In the last section, we prove the following theorem, regarding monotonicity of $|\lambda_n|$.
\begin{theorem}\label{ineq_eq}
For each $a\in\C^{m-1}$ there exists $M>0$ such that $|\lambda_n|<|\lambda_{n+1}|$ if $n\geq M$. 
\end{theorem}
This is a consequence of \eqref{asym_eq1}.

Finally, when $H$ is $\mathcal{PT}$-symmetric (i.e., $a\in\R^{m-1}$),  $u(z,\lambda)$ is an eigenfunction associated with an eigenvalue $\lambda$ if and only if $\overline{u(-\overline{z},\lambda)}$ is an eigenfunction associated with the eigenvalue $\overline{\lambda}$. Thus, the eigenvalues either appear in complex conjugate pairs, or else are real. So Theorem \ref{ineq_eq} implies the following.
\begin{theorem}\label{main_theorem}
Suppose that $a\in\R^{m-1}$. Then 
the eigenvalues $\lambda$ of $H$ are all real and positive, with only finitely many exceptions.
\end{theorem}

For the rest of the Introduction, we will mention a brief history of  problem \eqref{ptsym}. 

In recent years, these $\mathcal{PT}$-symmetric operators have gathered considerable attention, because ample numerical and asymptotic studies suggest that many of such operators have real eigenvalues only even though they are not self-adjoint. In particular, the differential operators $H$ with some polynomial potential $V$ and with the boundary condition \eqref{bdcond} have been considered by Bessis and Zinn-Justin (not in print), Bender and Boettcher \cite{Bender} and many other physicists \cite{Bender2,Bender5,CGM,Dorey2,mez,Ali1,Shin1,Simon,Znojil}. 

Around 1992 Bessis and Zinn-Justin (not in print) conjectured that when $V(z)=iz^3+\beta z^2$, $\beta\in\R$, the eigenvalues are all real and positive, and in 1998, Bender and Boettcher \cite{Bender} conjectured that when $V(z)=(iz)^m+\beta z^2$, $\beta\in\R$, the eigenvalues are all real and positive. Many numerical, asymptotic and analytic studies support these conjectures (see, e.\ g.,  \cite{Bender2,Bender5,CGM,Dorey2,mez,Ali1,Shin1,Simon,Znojil} and references therein and below).

The first rigorous proof of reality and positivity of the eigenvalues of some non-self-adjoint $H$ in \eqref{ptsym} was given by Dorey, Dunning and Tateo \cite{Dorey} in 2001. They  proved that the eigenvalues of $H$ with the potential $V(z)=-(iz)^{2m}-\alpha(iz)^{m-1}+\frac{\ell(\ell+1)}{z^2}$, $m,\,\alpha, \ell\in\R$, are all real if $m>1$ and $\alpha<m+1+|2\ell+1|$, and positive if $m>1$ and $\alpha<m+1-|2\ell+1|$. 

Then in 2002 the present author \cite{Shin} extended the polynomial potential results of Dorey, Dunning and Tateo to more general polynomial cases, by adapting the method in \cite{Dorey}. Namely, when $V(z)=-(iz)^m-P_{m-1}(iz)$, the eigenvalues are all real and positive, provided that for some $1\leq j\leq \frac{m}{2}$ the coefficients of the real polynomial $P_{m-1}$ satisfy $(j-k)a_k\geq 0$ for all $1\leq k\leq m-1$. 

However, there are some $\mathcal{PT}$-symmetric polynomial potentials that produce non-real eigenvalues. 
 Delabaere and Pham \cite{Pham}, and Delabaere and Trinh \cite{Delabaere} studied  the potential $iz^3+\gamma iz$ and showed that a pair of  
non-real eigenvalues develops for large negative $\gamma$. Moreover,  Handy \cite{Handy2}, and Handy, Khan, Wang and Tymczak \cite{Handy1} showed that the same potential admits a pair of non-real eigenvalues for small negative values of $\gamma\approx -3.0$. Also, Bender, Berry, Meisinger, Savage and Simsek \cite{Bender-1} considered the problem with the potential $V(z)=z^4+iAz$, $A\in\R$, under decaying boundary conditions at both ends of the real axis, and their numerical study showed that more and more non-real eigenvalues develop as $|A|\to\infty$. So without any restrictions on the coefficients $a_k$, Theorem \ref{main_theorem} is the most general result one can expect about reality of eigenvalues.

Also, the method used to prove  Theorem \ref{main_theorem} in this paper is new. The method used in \cite{Dorey, Shin} is useful in proving reality of {\it all} eigenvalues, but I think that  some critical arguments in proving reality of eigenvalues in  \cite{Dorey, Shin} cannot be applied to the cases when some non-real eigenvalues exist. The asymptotic expansion \eqref{asym_eq1} itself is interesting, and also \eqref{asym_eq1} implies Theorem \ref{ineq_eq}. Note that \eqref{Bender_exp} is not enough to conclude Theorem \ref{ineq_eq}. Finally,  Theorem \ref{ineq_eq} and $\mathcal{PT}$-symmetry of $H$ explained right before Theorem \ref{main_theorem} above imply the partial reality of the eigenvalues in  Theorem \ref{main_theorem}.

\section{Properties of the solutions}
\label{prop_sect}
In this section, we introduce work of Hille \cite{Hille} and Sibuya \cite{Sibuya} about properties of the solutions of \eqref{ptsym}. 

First, we scale equation \eqref{ptsym} because many facts that we need later are stated for the scaled equation. Let $u$ be a solution of (\ref{ptsym}) and let $v(z,\lambda)=u(-iz,\lambda)$. Then $v$ solves
\begin{equation}\label{rotated}
-v^\dd(z,\lambda)+[z^m+P_{m-1}(z)+\lambda]v(z,\lambda)=0,
\end{equation}
where $m\geq 3$ and $P_{m-1}$  is a  polynomial (possibly, $P_{m-1}\equiv 0$) of the form \eqref{poly_form}.

\begin{figure}[t]
    \begin{center}
    \includegraphics[width=.4\textwidth]{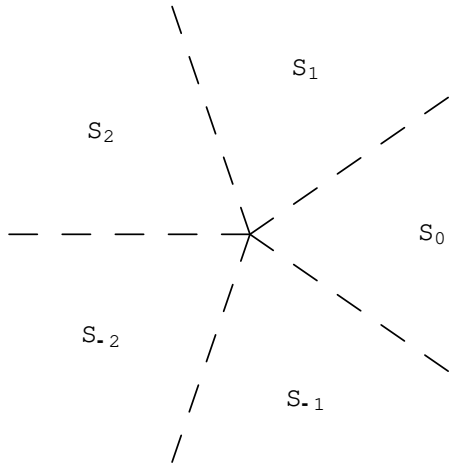}
    \end{center}
 \vspace{-.5cm}
\caption{The Stokes sectors for $m=3$. The dashed rays represent $\arg z=\pm\frac{\pi}{5},\,\pm\frac{3\pi}{5},\, \pi.$}\label{f:graph1}
\end{figure}

Since we scaled the argument of $u$, we must rotate the boundary conditions. We state them in a more general context by using the following definition.
\begin{definition}
{\rm {\it The Stokes sectors} $S_k$ of the equation (\ref{rotated})
are
$$ S_k=\left\{z\in \C:\left|\arg (z)-\frac{2k\pi}{m+2}\right|<\frac{\pi}{m+2}\right\}\quad\text{for}\quad k\in \Z.$$ }
\end{definition}
See Figure \ref{f:graph1}.
It is known from Hille \cite[\S 7.4]{Hille} that every nonconstant solution of (\ref{rotated}) either decays to zero or blows up exponentially, in each Stokes sector $S_k$. 
 That is, one has the following result.

\begin{lemma}[\protect{\cite[\S 7.4]{Hille}}]\label{gen_pro}
${} $

\begin{itemize}
\item [(i)] For each $k\in\Z$, every solution $v$ of (\ref{rotated}) (with no boundary conditions imposed) is asymptotic to 
\begin{equation} \label{asymp-formula}
(const.)z^{-\frac{m}{4}}\exp\left[\pm \int^z \left[\xi^m+P_{m-1}(\xi)+\lambda\right]^{\frac{1}{2}}\,d\xi\right] 
\end{equation}
as $z \rightarrow \infty$ in every closed subsector of $S_k$.

\item [(ii)] If a nonconstant solution $v$ of \eqref{rotated} decays in $S_k$, it must blow up in $S_{k-1}\cup S_{k+1}$. However, when $v$ blows up in $S_k$, $v$ need not be decaying in $S_{k-1}$ or in $S_{k+1}$.
\end{itemize}
\end{lemma}
Lemma \ref{gen_pro} (i) implies that if $v$ decays along one ray in $S_k$, then it decays along all rays in $S_k$. Also, if $v$ blows up along one ray in $S_k$, then it blows up along all rays in $S_k$.
Thus, since the rotation $z\mapsto iz$ maps the two rays in \eqref{bdcond} onto the center rays of $S_{-1}$ and  $S_1$, 
$$\text{the boundary conditions on $u$ in \eqref{ptsym} mean that $v$ decays in $S_{-1}\cup S_1$}.$$

Next we will introduce Sibuya's results, but first we define  a sequence of complex numbers $b_j$ in terms of the $a_k$ and $\lambda$, as follows. For $\lambda \in \C$ fixed, we expand
\begin{eqnarray}
&&(1+a_1z^{-1}+a_2z^{-2}+\cdots+a_{m-1}z^{1-m}+\lambda z^{-m})^{1/2}\nonumber\\
&=&1+\sum_{k=1}^{\infty}{\frac{1}{2}\choose{k}}\left( a_1z^{-1}+a_2z^{-2}+\cdots+a_{m-1}z^{1-m}+\lambda z^{-m}\right)^k\nonumber\\
&=&1+\sum_{j=1}^{\infty}\frac{b_j(a,\lambda)}{z^j},\qquad\text{for large}\quad|z|.\label{b_def}
\end{eqnarray}
Note that $b_1,\,b_2,\,\ldots,\,b_{m-1}$ do not depend on $\lambda$, so we write $b_j(a)=b_j(a,\lambda)$ for $j=1,\,2,\dots,\,m-1$. So the above expansion without the $\lambda z^{-m}$ term still gives $b_j$ for $1\leq j\leq m-1$. We further define $r_m=-\frac{m}{4}$ if $m$ is odd, and $r_m=-\frac{m}{4}-b_{\frac{m}{2}+1}(a)$ if $m$ is even.

 The following theorem is a special case of Theorems 6.1, 7.2, 19.1 and 20.1 of Sibuya \cite{Sibuya} that is the main ingredient of the proofs of the main results in this paper. 
\begin{theorem}\label{prop}
Equation (\ref{rotated}), with $a\in \C^{m-1}$, admits a solution  $f(z,a,\lambda)$ with the following properties.
\begin{enumerate}
\item[(i)] $f(z,a,\lambda)$ is an entire function of $z,a $ and $\lambda$.
\item[(ii)] $f(z,a,\lambda)$ and $f^\d(z,a,\lambda)=\frac{\partial}{\partial z}f(z,a,\lambda)$ admit the following asymptotic expansions. Let $\varepsilon>0$. Then
\begin{align}
f(z,a,\lambda)=&\qquad z^{r_m}(1+O(z^{-1/2}))\exp\left[-F(z,a,\lambda)\right],\nonumber\\
f^\d(z,a,\lambda)=&-z^{r_m+\frac{m}{2}}(1+O(z^{-1/2}))\exp\left[-F(z,a,\lambda) \right],\nonumber
\end{align}
as $z$ tends to infinity in  the sector $|\arg z|\leq \frac{3\pi}{m+2}-\varepsilon$, uniformly on each compact set of $(a,\lambda)$-values . 
Here
\begin{equation}\nonumber
F(z,a,\lambda)=\frac{2}{m+2}z^{\frac{m}{2}+1}+\sum_{1\leq j<\frac{m}{2}+1}\frac{2}{m+2-2j}b_j(a) z^{\frac{1}{2}(m+2-2j)}.
\end{equation}
\item[(iii)] Properties \textup{(i)} and \textup{(ii)} uniquely determine the solution $f(z,a,\lambda)$ of (\ref{rotated}).
\item[(iv)] For each fixed $a\in\C^{m-1}$ and $\delta>0$, $f$ and $f^\d$ also admit the asymptotic expansions, 
\begin{align}
f(0,a,\lambda)=&[1+o(1)]\lambda^{-1/4}\exp\left[L(a,\lambda)\right],\label{eq1}\\
f^\d(0,a,\lambda)=&-[1+o(1)]\lambda^{1/4}\exp\left[L(a,\lambda)\right],\label{eq2}
\end{align}
as $\lambda\to\infty$ in the sector $|\arg(\lambda)|\leq\pi-\delta$, where 
\begin{align}
L(a,\lambda)=\left\{\begin{array}{rl}
                    \int_0^{+\infty}\left(\sqrt{t^m+P_{m-1}(t)+\lambda}- t^{\frac{m}{2}}-\sum_{j=1}^{\frac{m+1}{2}}b_j(a)t^{\frac{m}{2}-j}\right)dt&\text{ if $m$ is odd,}\\
                   \int_0^{+\infty}\left(\sqrt{t^m+P_{m-1}(t)+\lambda}- t^{\frac{m}{2}}-\sum_{j=1}^{\frac{m}{2}}b_j(a)t^{\frac{m}{2}-j}-\frac{b_{\frac{m}{2}+1}(a)}{t+1}\right)dt&\text{ if $m$ is even.}
                    \end{array}
                         \right. \nonumber
\end{align}
\item[(v)] The entire functions  $\lambda\mapsto f(0,a,\lambda)$ and $\lambda\mapsto f^\d(0,a,\lambda)$ have orders $\frac{1}{2}+\frac{1}{m}$.
\end{enumerate}
\end{theorem}
\begin{proof}
In  Sibuya's book \cite{Sibuya}, see Theorem 6.1 for a proof of (i) and (ii); Theorem 7.2 for a proof of (iii); and Theorem 19.1 for a proof of (iv).  Moreover, (v) is a consequence of (iv) along with Theorem 20.1. Note that properties (i), (ii) and (iv) are summarized on pages 112--113 of Sibuya \cite{Sibuya}. 
\end{proof}

Using this theorem, Sibuya \cite[Theorem 19.1]{Sibuya} also showed the following corollary that will be  useful  later on.
\begin{corollary}\label{lemma_decay}
Let $a\in\C^{m-1}$ be fixed. Then
$L(a,\lambda)=K_{m}\lambda^{\frac{1}{2}+\frac{1}{m}}(1+o(1))$ as $\lambda$ tends to infinity in the sector $|\arg \lambda|\leq \pi-\delta$, and hence
\begin{equation}\label{re_part}
\Re \left(L(a,\lambda)\right)=K_{m}\cos\left(\frac{m+2}{2m}\arg(\lambda)\right)|\lambda|^{\frac{1}{2}+\frac{1}{m}}(1+o(1))
\end{equation}
as $\lambda\to\infty$ in the sector $|\arg (\lambda)|\leq \pi-\delta$.

In particular, $\Re \left(L(a,\lambda)\right)\to+\infty$ as $\lambda\to\infty$ in any closed subsector of the sector $|\arg(\lambda)|<\frac{m\pi}{m+2}$. In addition, $\Re \left(L(a,\lambda)\right)\to-\infty$ as $\lambda\to\infty$ in any closed subsector of the sectors $\frac{m\pi}{m+2}<|\arg(\lambda)|<\pi-\delta$.
\end{corollary}
\begin{proof}
This asymptotic expansion will be clear from Lemma \ref{asy_lemma} below, or alternatively, see \cite[Theorem 19.1]{Sibuya} for a proof.
\end{proof}

Based on the above Corollary,
Sibuya \cite[Theorem 29.1]{Sibuya} also proved the following asymptotic expansion of the eigenvalues.
\begin{equation}\label{ak}
\lambda_n= \omega^{m}\left(\frac{(-2n+1)\pi}{2K_{m} \sin \left(\frac{2\pi}{m}\right)} \right)^{\frac{2m}{m+2}}[1+o(1)],\quad\text{as}\quad n\rightarrow \infty,
\end{equation}
where
$$\omega=\exp\left[\frac{2\pi i}{m+2}\right].$$
Notice that in this paper we consider the boundary conditions of the scaled equation (\ref{rotated}) where $v$ decays in $S_{-1}\cup S_1$, while 
  Sibuya studies  equation (\ref{rotated}) with  boundary conditions such that $v$ decays in $S_0\cup S_2$. The factor $\omega^{m}$ in our formula (\ref{ak}) is due to this scaling of the problem. 

\begin{remark}
{\rm Throughout this paper, we will deal with numbers like $\left(\omega^{\nu}\lambda\right)^{s}$ for some $s\in\R$, and $\nu\in\C$. As usual, we will use
$$\omega^{\nu}=\exp\left[\nu \frac{2\pi i}{m+2}\right]$$
and if $\arg(\lambda)$ is specified, then
$$\arg\left(\left(\omega^{\nu}\lambda\right)^{s}\right)=s\left[\arg(\omega^{\nu})+\arg(\lambda)\right]=s\left[\Re(\nu)\frac{2\pi}{m+2}+\arg(\lambda)\right],\quad s\in\R.$$
If $s\not\in\Z$ then the branch of $\lambda^s$ is chosen to be the negative real axis.
}
\end{remark}

Next, we provide an improved asymptotic expansion of $L$. We will use this new asymptotic expansion of $L$ to improve the asymptotic expansion \eqref{ak} of the eigenvalues.   
\begin{lemma}\label{asy_lemma}
Let $m\geq 3$ and $a\in\C^{m-1}$ be fixed. Then there exist constants $K_{m,j}(a)\in\C$, $0\leq j \leq\frac{m}{2}+1$, such that 
\begin{equation}
L(a,\lambda)=\left\{\begin{array}{rl}
&\sum_{j=0}^{\frac{m+1}{2}}K_{m,j}(a)\lambda^{\frac{1}{2}+\frac{1-j}{m}}+O\left(|\lambda|^{-\frac{1}{2m}}\right)\,\,\text{if $m$ is odd,}\\
&\sum_{j=0}^{\frac{m}{2}+1}K_{m,j}(a)\lambda^{\frac{1}{2}+\frac{1-j}{m}}-\frac{b_{\frac{m}{2}+1}(a)}{m}\ln(\lambda)+O\left(|\lambda|^{-\frac{1}{m}}\right)\,\,\text{if $m$ is even,}
\end{array}
\right.\nonumber
\end{equation}
as $\lambda\to\infty$ in the sector $|\arg(\lambda)|\leq \pi-\delta$. 
\end{lemma}
\begin{proof}
The function $L(a,\lambda)$ is defined as an integral over $0\leq t<+\infty$ in Theorem \ref{prop}. We will rotate the contour of integration  using  Cauchy's integral formula. In doing so, we need to justify that the integrand in the definition of $L(a,\lambda)$ is analytic in some domain in the complex plane.
 
Let $0<\delta<\frac{\pi}{m+2}$ be a fixed number. Suppose that $0\leq\arg(\lambda)\leq \pi-\delta$. Then
if $0\leq\arg(t)\leq\frac{1}{m}\arg(\lambda)$,  there exists $M_0>0$ such that 
\begin{equation}\nonumber
-\pi<-\frac{\delta}{2}\leq\arg(t^m+P_{m-1}(t))\leq \arg(\lambda)+\frac{\delta}{2}\leq\pi-\frac{\delta}{2},
\end{equation}
provided that $|t|\geq M_0$. Since $t^m+P_{m-1}(t)$ lies in a large disk centered at the origin for $|t|\leq M_0$, we see that for all $\lambda$ with $|\lambda|$ large, we have that $-\frac{\delta}{2}<\arg(t^m+P_{m-1}(t)+\lambda)<\pi-\frac{\delta}{2}$ and $|t^m+P_{m-1}(t)+\lambda|>0$ for all $t$ in the sector $0\leq\arg(t)\leq\frac{1}{m}\arg(\lambda)$, and hence $ \sqrt{t^m+P_{m-1}(t)+\lambda}$ is analytic in the sector $0\leq\arg(t)\leq\frac{1}{m}\arg(\lambda)$ if $\lambda$ lies outside a large disk and in the sector $0\leq\arg(\lambda)\leq \pi-\delta$.

Let 
\begin{align}
Q(t,a,\lambda)=\left\{
                    \begin{array}{rl}
                    &\sqrt{t^m+P_{m-1}(t)+\lambda}- t^{\frac{m}{2}}-\sum_{j=1}^{\frac{m+1}{2}}b_j(a)t^{\frac{m}{2}-j} \quad \text{if $m$ is odd,}\\
                   &\sqrt{t^m+P_{m-1}(t)+\lambda}- t^{\frac{m}{2}}-\sum_{j=1}^{\frac{m}{2}}b_j(a)t^{\frac{m}{2}-j}-\frac{b_{\frac{m}{2}+1}(a)}{t+1}  \quad \text{if $m$ is even.}
                    \end{array}
                         \right. \nonumber
\end{align}
 Then, since $|Q(t,a,\lambda)|=O\left(|t|^{-\frac{m}{2}}\right)$ as $t$ tends to infinity in the sector $0\leq\arg(t)\leq\frac{1}{m}\arg(\lambda)$,  we have by Cauchy's integral formula, upon substituting $t=\lambda^{\frac{1}{m}}\tau$ for all $\lambda$ with $|\lambda|$ large enough, 
\begin{equation}\label{asy_eql}
L(a,\lambda)=\int_{0}^{+\infty}Q(t,a,\lambda)\, dt
=\lambda^{\frac{1}{m}}\int_0^{+\infty}Q(\lambda^{\frac{1}{m}}\tau,a,\lambda)\, d\tau,
\end{equation}
where
\begin{align}
Q(\lambda^{\frac{1}{m}}& \tau, a,\lambda)\nonumber\\
=&\left\{
                    \begin{array}{rl}
                    &\lambda^{\frac{1}{2}}\left(\sqrt{\tau^m+1+\frac{P_{m-1}(\lambda^{\frac{1}{m}}\tau)}{\lambda}}- \tau^{\frac{m}{2}}-\sum_{j=1}^{\frac{m+1}{2}}b_j(a)\frac{\tau^{\frac{m}{2}-j}}{\lambda^{\frac{j}{m}}}\right)\,\, \text{if $m$ is odd,}\\
                   &\lambda^{\frac{1}{2}}\left(\sqrt{\tau^m+1+\frac{P_{m-1}(\lambda^{\frac{1}{m}}\tau)}{\lambda}}- \tau^{\frac{m}{2}}-\sum_{j=1}^{\frac{m}{2}}b_j(a)\frac{\tau^{\frac{m}{2}-j}}{\lambda^{\frac{j}{m}}}-\frac{\lambda^{-\frac{1}{2}}b_{\frac{m}{2}+1}(a)}{\lambda^{\frac{1}{m}}\tau+1}\right)
 \,\, \text{if $m$ is even.}
                    \end{array}
                         \right. \nonumber
\end{align}

Similarly, \eqref{asy_eql} holds for  $-\pi+\delta\leq\arg(\lambda)\leq 0$.

Next, we examine the following square root in $Q(\lambda^{\frac{1}{m}}\tau, a,\lambda)$:
\begin{align}
\sqrt{\tau^m+1+\frac{P_{m-1}(\lambda^{\frac{1}{m}}\tau)}{\lambda}}
&=\sqrt{\tau^m+1}\sqrt{1+\frac{P_{m-1}(\lambda^{\frac{1}{m}}\tau)}{\lambda(\tau^m+1)}}\nonumber\\
&=\sqrt{\tau^m+1}\left(1+\sum_{k=1}^{\infty}{{\frac{1}{2}}\choose{k}}\left(\frac{P_{m-1}(\lambda^{\frac{1}{m}}\tau)}{\lambda(\tau^m+1)}\right)^k\right)\nonumber\\
&\overset{let}{=}\sqrt{\tau^m+1}+\sum_{j=1}^{\infty}\frac{g_j(\tau)}{\lambda^{\frac{j}{m}}},\nonumber
\end{align}
where $g_j(\tau)$ are functions such that $g_j(\tau)$ are all integrable on $[0,\, R]$ for any $R>0$. Moreover, by the definition of $b_j$ in \eqref{b_def}, we see that for $1\leq j\leq m-1$,
$$g_j(\tau)=\sum_{k=1}^{j}\frac{b_{j,k}(a)\tau^{mk-j}}{\left(\tau^m+1\right)^{k-\frac{1}{2}}}\quad\text{for some constants $b_{j,k}(a)$ such that $\sum_{k=1}^jb_{j,k}(a)=b_j(a)$}.$$
Thus,
\begin{align}
g_j(\tau)-b_j(a)\tau^{\frac{m}{2}-j}&=\sum_{k=1}^jb_{j,k}(a)\left(\frac{\tau^{mk-j}}{\left(\tau^m+1\right)^{k-\frac{1}{2}}}-\tau^{\frac{m}{2}-j}\right)\nonumber\\
&\underset{\tau\to\infty}{=}\sum_{k=1}^jb_{j,k}(a)\tau^{\frac{m}{2}-j}O\left(\frac{1}{\tau^m}\right)\nonumber\\
&\underset{\tau\to\infty}{=}O\left(\frac{1}{\tau^{\frac{m}{2}+j}}\right)\quad\text{for all $1\leq j\leq\frac{m+1}{2}$}.\nonumber
\end{align}

So $\int_0^{\infty}\left|g_j(\tau)-b_j(a)\tau^{\frac{m}{2}-j}\right|\,d\tau<+\infty$ for all $1\leq j\leq\frac{m+1}{2}$.
Next, when $m$ is even and $j=\frac{m}{2}+1$, 
we write
\begin{align}
&\int_0^{\infty}\left(g_{\frac{m}{2}+1}(\tau)-\frac{b_{\frac{m}{2}+1}(a)}{\tau+\lambda^{-\frac{1}{m}}}\right)\, d\tau\nonumber\\
&=\int_0^{\infty}\left(g_{\frac{m}{2}+1}(\tau)-\frac{b_{\frac{m}{2}+1}(a)}{\tau+1}\right)\, d\tau+b_{\frac{m}{2}+1}(a)\int_0^{\infty}\left(\frac{1}{\tau+1}-\frac{1}{\tau+\lambda^{-\frac{1}{m}}}\right)\, d\tau\nonumber\\
&\overset{let}{=}K_{m,\frac{m}{2}+1}(a)-\frac{b_{\frac{m}{2}+1}(a)}{m}\ln(\lambda),\nonumber
\end{align}
where we take $\Im(\ln(\lambda))=\arg(\lambda)\in(-\pi,\pi)$.

Thus, we have that
\begin{equation}
L(a,\lambda)=\left\{\begin{array}{rl}
&\sum_{j=0}^{\frac{m+1}{2}}K_{m,j}(a)\lambda^{\frac{1}{2}+\frac{1-j}{m}}+O\left(|\lambda|^{-\frac{1}{2m}}\right)\,\,\text{if $m$ is odd,}\\
&\sum_{j=0}^{\frac{m}{2}+1}K_{m,j}(a)\lambda^{\frac{1}{2}+\frac{1-j}{m}}-\frac{b_{\frac{m}{2}+1}(a)}{m}\ln(\lambda)+O\left(|\lambda|^{-\frac{1}{m}}\right)\,\,\text{if $m$ is even,}
\end{array}
\right.\nonumber
\end{equation}
as $\lambda\to\infty$ in the sector $|\arg(\lambda)|\leq\pi-\delta$, where
\begin{align}
 K_{m,0}(a)&=K_{m}=\int_0^{\infty}\left(\sqrt{1+t^m}-\sqrt{t^m}\right)\, dt>0\quad\text{for all $m\geq 3$},\nonumber\\
K_{m,j}(a)&=\int_0^{\infty}\left(g_j(t)-b_j(a)t^{\frac{m}{2}-j}\right)\,dt\quad\text{for all $1\leq j\leq\frac{m+1}{2}$},\label{def_K}\\
K_{m,\frac{m}{2}+1}(a)&=\int_0^{\infty}\left(g_{\frac{m}{2}+1}(t)-\frac{b_{\frac{m}{2}+1}(a)}{t+1}\right)\,dt\quad\text{when $m$ is even}.\nonumber
\end{align}
This completes the proof.
\end{proof}

\section{Eigenvalues are zeros of an entire function}\label{entire_sect}

In this section, we will prove that the eigenvalues are zeros of an entire function.
  
First, we let
\begin{equation}\nonumber
G^k(a):=(\omega^{-k}a_1, \omega^{-2k}a_2,\ldots,\omega^{-(m-1)k}a_{m-1})\quad \text{for}\quad k\in \Z.
\end{equation}
Then recall that the function
 $f(z,a,\lambda)$  in Theorem \ref{prop} solves (\ref{rotated}) and decays to zero exponentially as $z\rightarrow \infty$ in $S_0$, and  blows up in $S_{-1}\cup S_1$. Next, one can check that the function
$$f_k(z,a,\lambda):=f(\omega^{-k}z,G^k(a),\omega^{-mk}\lambda),$$
 which is obtained by scaling $f(z,G^k(a),\omega^{-mk}\lambda)$ in the $z$-variable, also solves (\ref{rotated}). It is clear that $f_0(z,a,\lambda)=f(z,a,\lambda)$. Also,  $f_k(z,a,\lambda)$ decays in $S_k$ and blows up in $S_{k-1}\cup S_{k+1}$ since $f(z,G^k(a),\omega^{-mk}\lambda)$ decays in $S_0$. Since no nonconstant solution decays in two consecutive Stokes sectors (see Lemma \ref{gen_pro} (ii)), $f_{k}$ and $f_{k+1}$ are linearly independent and hence any solution of (\ref{rotated}) can be expressed as a linear combination of these two. Especially,  there exist some coefficients $C(a,\lambda)$ and $\widetilde{C}(a,\lambda)$ such that
\begin{equation}\label{stokes}
f_{-1}(z,a,\lambda)=C(a,\lambda)f_0(z,a,\lambda)+\widetilde{C}(a,\lambda)f_{1}(z,a,\lambda).
\end{equation}

We then see that 
\begin{equation}\label{C_def}
C(a,\lambda)=\frac{W_{-1,1}(a,\lambda)}{W_{0,1}(a,\lambda)}\quad\text{and}\quad \widetilde{C}(a,\lambda)=-\frac{W_{-1,0}(a,\lambda)}{W_{0,1}(a,\lambda)},
\end{equation}
where $W_{j,k}=f_jf_k^\d -f_j^\d f_k$ is the Wronskian of $f_j$ and $f_k$. Since both $f_j$ and $f_k$ are solutions of the same linear equation (\ref{rotated}), we know that the Wronskians are constant functions of $z$. Also, since  $f_k$ and $f_{k+1}$ are linearly independent,  $W_{k,k+1}\not=0$ for all $k\in \Z$. Moreover, we have the following lemma that is useful  later on.
\begin{lemma}
Suppose $k,\,j\in\Z$. Then
\begin{equation}\label{kplus1}
W_{k+1,j+1}(a,\lambda)=\omega^{-1}W_{k,j}(G(a),\omega^2\lambda),
\end{equation}
and $W_{0,1}(a,\lambda)=2\omega^{\mu(a)}$, where
\begin{eqnarray}
\mu(a)=\left\{
              \begin{array}{rl}
              \frac{m}{4}  \quad &\text{if $m$ is odd,}\\
              \frac{m}{4}- b_{\frac{m}{2}+1}(a) \quad &\text{if $m$ is even.}
              \end{array}
                         \right. \nonumber
\end{eqnarray}
Moreover, 
\begin{equation}\nonumber
\widetilde{C}(a,\lambda)=-\frac{W_{-1,0}(a,\lambda)}{W_{0,1}(a,\lambda)}
=-\omega \frac{W_{0,1}(G^{-1}(a),\omega^{-2}\lambda)}{W_{0,1}(a,\lambda)}=-\omega^{1+2\nu(a)},
\end{equation}
where
\begin{eqnarray}
\nu(a)=\left\{
                         \begin{array}{rl}
                         0 & \quad \text{if $m$ is odd,}\\
                         b_{\frac{m}{2}+1}(a)& \quad \text{if $m$ is even.}
                          \end{array}
                         \right. \label{def_nu}
\end{eqnarray}
\end{lemma}
\begin{proof}
See Sibuya \cite[pages 116-118]{Sibuya} for proof. Here, we mention that by \eqref{b_def}, we have $b_{\frac{m}{2}+1}(G^{-1}(a))=-b_{\frac{m}{2}+1}(a)$ and hence $\nu(G^{-1}(a))=-\nu(a)$.
\end{proof}

Now we can identify the eigenvalues of $H$ as the zeros of 
the entire function $\lambda\mapsto C(a,\lambda)$.
\begin{theorem}\label{iden_thm}
For each fixed $a\in\C^{m-1}$, the function $\lambda\mapsto C(a,\lambda)$ is  entire.
Moreover, $\lambda$ is an eigenvalue of $H$ if and only if $C(a,\lambda)=0$.
\end{theorem}
\begin{proof}
Since $W_{0,1}(a,\lambda)\not=0$ and since $W_{-1, 1}(a,\lambda)$ is a Wronskian of two entire functions, it is clear from \eqref{C_def} that  $C(a,\lambda)$ is an entire function of $\lambda$ for each fixed $a\in\C^{m-1}$.

Next, suppose that $\lambda$ is an eigenvalue of $H$ with a corresponding eigenfunction $u$, then  the scaled eigenfunction $v(z,\lambda)=u(-iz,\lambda)$ solves \eqref{rotated} and decays in $S_{-1}\cup S_1$. Hence, $v$ is a (nonzero) constant multiple of $f_1$ since both decays in $S_1$. Similarly, $v$ is also a constant multiple of $f_{-1}$. Thus, $f_{-1}$ is a constant multiple of $f_1$, implying $C(a,\lambda)=0$.

Conversely, if $C(a,\lambda)=0$, then $f_{-1}$ is a constant multiple of $f_1$, and hence $f_1$ also decays in $S_{-1}$. Thus, $f_1$ decays in $S_{-1}\cup S_{ 1}$ and is  a scaled eigenfunction with the eigenvalue $\lambda$.
\end{proof}

Moreover, the following is an easy consequence of \eqref{stokes}: For each $k\in\Z$ we have 
\begin{equation}\label{kdef}
W_{-1,k}(a,\lambda)=C(a,\lambda)W_{0,k}(a,\lambda)+\widetilde{C}(a)W_{1,k}(a,\lambda),
\end{equation}
where we use $\widetilde{C}(a)$ for $\widetilde{C}(a,\lambda)$ since it is independent of $\lambda$.

\section{Asymptotic expansions of $C(a,\lambda)$}\label{asymp_sect}
In this section, we  provide asymptotic expansions of the entire function $C(a,\lambda)$ as $\lambda\to\infty$ along all possible rays to infinity in the complex plane.

First, we provide an asymptotic expansion of  the Wronskian of $f_0$ and $f_j$ in preparation for providing an asymptotic expansion of  $C(a,\lambda)$.
\begin{lemma}\label{lemma7}
Suppose that $1\leq j\leq \frac{m}{2}+1$. Then for each $a\in\C^{m-1}$,
\begin{equation}\label{sec_eq1}
W_{0,j}(a,\lambda)=[2i\omega^{-\frac{j}{2}}+o(1)]\exp\left[L(G^{j}(a),\omega^{2j-m-2}\lambda)+L(a,\lambda)\right],
\end{equation}
as $\lambda\to\infty$ in the sector
\begin{equation}\label{sector0}
-\pi+\delta\leq \pi-\frac{4j\pi}{m+2}+\delta \leq \arg(\lambda)\leq \pi-\delta.
\end{equation}
\end{lemma}
\begin{proof}
We fix $1\leq j\leq \frac{m}{2}+1$. Then,
\begin{align}
W_{0,j}(a,\lambda)
&\,\,\,=f_0(z,a,\lambda)f_j^\d(z,a,\lambda)-f_0^\d(z,a,\lambda)f_j(z,a,\lambda)\nonumber\\
&\,\,\,=\omega^{-j}f(0,a,\lambda)f^\d(0,G^{j}(a),\omega^{2j-m-2}\lambda)-f^\d(0,a,\lambda)f(0,G^{j}(a),\omega^{2j-m-2}\lambda)\nonumber\\
&\underset{\lambda\to\infty}{=}-\left[\omega^{-j}\omega^{\frac{2j-m-2}{4}}-\omega^{-\frac{2j-m-2}{4}}+o(1)\right]\exp\left[L(G^{j}(a),\omega^{2j-m-2}\lambda)+L(a,\lambda)\right]\nonumber\\
&\underset{\lambda\to\infty}{=}[2i\omega^{-\frac{j}{2}}+o(1)]\exp\left[L(G^{j}(a),\omega^{2j-m-2}\lambda)+L(a,\lambda)\right],\nonumber
\end{align}
where we used \eqref{eq1} and \eqref{eq2} with
\begin{equation}\nonumber
|\arg (\lambda)|\leq \pi-\delta\quad\text{and}\quad |\arg (\omega^{2j-m-2}\lambda)|\leq \pi-\delta,
\end{equation}
which is, \eqref{sector0}. Here we also used $j\leq\frac{m}{2}+1$.
\end{proof}

Next, we provide an asymptotic expansion of $W_{-1,1}(a,\lambda)$ as $\lambda\to\infty$ along the rays near the negative real axis. Notice from \eqref{C_def} that $W_{-1,1}(a,\lambda)=W_{0,1}(a,\lambda)C(a,\lambda)$. Also, $W_{0,1}(a,\lambda)$ is a nonzero constant function of  $\lambda$. So from these one gets an asymptotic expansion of $C(a,\lambda)$.

\begin{theorem}
For each fixed $a\in\C^{m-1}$ and $0<\delta<\frac{\pi}{m+2}$,
\begin{equation}\label{asy_1}
W_{-1,1}(a,\lambda)=[2i+o(1)]\exp\left[L(G^{-1}(a),\omega^{-2}\lambda)+L(G(a),\omega^{-m}\lambda)\right],
\end{equation}
as $\lambda\to \infty$ along the rays in the sector
\begin{equation}\label{sector1}
\pi-\frac{4\pi}{m+2}+\delta\leq \arg(\lambda)\leq \pi+\frac{4\pi}{m+2}-\delta.
\end{equation}
Moreover, there exists a constant $M_1>0$ such that $W_{-1,1}(a,\lambda)\not=0$ for all $\lambda$ in the sector \eqref{sector1} if $|\lambda|\geq M_1$.
\end{theorem}
\begin{proof}
This is an easy consequence of Lemma \ref{lemma7} and equation \eqref{kplus1}.

The last assertion of the theorem is a consequence of the asymptotic expansion \eqref{asy_1}.
\end{proof}

The asymptotic expansion of $C(a,\lambda)$ in a sector near the positive real axis is obtained in the following theorem. 
\begin{theorem}\label{zero_thm}
Suppose that $m\geq 4$. Then for each fixed $a\in\C^{m-1}$ and $0<\delta<\frac{\pi}{m+2},$
\begin{align}
C(a,\lambda)=&[\omega^{\frac{1}{2}}+o(1)]\exp\left[L(G^{-1}(a),\omega^{-2}\lambda)-L(a,\lambda)\right]\nonumber\\
&+[\omega^{\frac{1}{2}+2\nu(a)}+o(1)]\exp\left[L(G(a),\omega^{2}\lambda)-L(a,\lambda)\right],\nonumber
\end{align}
as $\lambda\to\infty$ in the sector 
\begin{equation}\label{sector4}
\pi-\frac{4\lfloor\frac{m}{2}\rfloor\pi}{m+2}+\delta \leq \arg(\lambda)\leq \pi-\frac{4\pi}{m+2}-\delta.
\end{equation} 
\end{theorem}

\begin{proof}
Suppose $2\leq k\leq\frac{m}{2}$. Then from \eqref{kplus1}, \eqref{kdef} and Lemma \ref{lemma7},
\begin{align}
C(a,\lambda)=&\frac{W_{-1,k}(a,\lambda)}{W_{0,k}(a,\lambda)}-\widetilde{C}(a)\frac{W_{1,k}(a,\lambda)}{W_{0,k}(a,\lambda)}\nonumber\\
=&\frac{\omega W_{0,k+1}(G^{-1}(a),\omega^{-2}\lambda)}{W_{0,k}(a,\lambda)}
-\widetilde{C}(a)\frac{\omega^{-1}W_{0,k-1}(G(a),\omega^{2}\lambda)}{W_{0,k}(a,\lambda)}\nonumber\\
=&[\omega^{\frac{1}{2}}+o(1)]\frac{\exp\left[L(G^{k}(a),\omega^{2k-m-2}\lambda)+L(G^{-1}(a),\omega^{-2}\lambda)\right]}{\exp\left[L(G^{k}(a),\omega^{2k-m-2}\lambda)+L(a,\lambda)\right]}\nonumber\\
&-[\omega^{-\frac{1}{2}}+o(1)]\widetilde{C}(a)\frac{\exp\left[L(G^{k}(a),\omega^{2k-m-2}\lambda)+L(G(a),\omega^{2}\lambda)\right]}{\exp\left[L(G^{k}(a),\omega^{2k-m-2}\lambda)+L(a,\lambda)\right]}\nonumber\\
=&[\omega^{\frac{1}{2}}+o(1)]\exp\left[L(G^{-1}(a),\omega^{-2}\lambda)-L(a,\lambda)\right]\nonumber\\
&+[\omega^{-\frac{1}{2}}+o(1)]\omega^{1+2\nu(a)}\exp\left[L(G(a),\omega^{2}\lambda)-L(a,\lambda)\right],\nonumber
\end{align}
as $\lambda\to\infty$ such that
\begin{align}
-\pi<\pi-\frac{4(k+1)\pi}{m+2}+\delta &\leq \arg(\omega^{-2}\lambda)\leq \pi-\delta,\nonumber\\
\pi-\frac{4k\pi}{m+2}+\delta &\leq \arg(\lambda)\leq \pi-\delta,\nonumber\\
\pi-\frac{4(k-1)\pi}{m+2}+\delta &\leq \arg(\omega^{2}\lambda)\leq \pi-\delta,\nonumber
\end{align}
that is,
\begin{equation}\nonumber
\pi-\frac{4k\pi}{m+2}+\delta \leq \arg(\lambda)\leq \pi-\frac{4\pi}{m+2}-\delta,
\end{equation}
provided that $2\leq k\leq \frac{m}{2}.$
So in order to complete the proof, we choose 
\begin{eqnarray}
k=\left\{
              \begin{array}{rl}
              \frac{m-1}{2}  \quad & \text{if $m$ is odd and $m\geq 5$,}\\
              \frac{m}{2} \quad & \text{if $m$ is even.}
              \end{array}
                         \right. \nonumber
\end{eqnarray}
\end{proof}

The sectors \eqref{sector1} and \eqref{sector4} do not cover the entire complex plane near infinity.
The next theorem covers a sector in the upper half plane, connecting the sectors \eqref{sector1} and \eqref{sector4} in the upper half plane.

\begin{theorem}\label{thm_sector2}
Suppose that $a\in\C^{m-1}$ and $0<\delta<\frac{\pi}{m+2}$.  If $m\geq 4$ then
\begin{align}
C(a,\lambda)=&\left[\omega^{\frac{1}{2}}+o(1)\right]\exp\left[L(G^{-1}(a),\omega^{-2}\lambda)-L(a,\lambda)\right]\nonumber\\
&-[i\omega^{1+\mu(a)+4\nu(a)}+o(1)]\exp\left[-L(G^2(a),\omega^{2-m}\lambda)-L(a,\lambda)\right],\label{asy_asy}
\end{align}
as $\lambda\to\infty$ in the sector
\begin{equation}\label{sector2}
\pi-\frac{8\pi}{m+2}+\delta \leq \arg(\lambda)\leq \pi-\delta.
\end{equation}
If $m=3$ then 
\begin{align}
C(a,\lambda)=&[-\omega^{-2}+o(1)]\exp\left[L(G^{4}(a),\omega^{-2}\lambda)-L(a,\lambda)\right]\nonumber\\
&-[i\omega^{\frac{7}{4}}+o(1)]\exp\left[-L(G^2(a),\omega^{-1}\lambda)-L(a,\lambda)\right],\nonumber
\end{align}
as $\lambda\to\infty$ in the sector 
\begin{equation}\label{sector2-1}
-\frac{\pi}{5}+\delta \leq \arg(\lambda)\leq \pi-\delta.
\end{equation}

Moreover, if $m\geq 6$ then there exists a constant $M_2>0$ such that $C(a,\lambda)\not=0$ for all $\lambda$  in the sector \eqref{sector2} if $|\lambda|\geq M_2$. 
\end{theorem}
\begin{proof}
Suppose that $m\geq 4$. Then from  \eqref{kplus1}, \eqref{kdef} and Lemma \ref{lemma7},
\begin{align}
C(a,\lambda)=&\frac{W_{-1,2}(a,\lambda)}{W_{0,2}(a,\lambda)}-\widetilde{C}(a)\frac{W_{1,2}(a,\lambda)}{W_{0,2}(a,\lambda)}\nonumber\\
=&\frac{\omega W_{0,3}(G^{-1}(a),\omega^{-2}\lambda)}{W_{0,2}(a,\lambda)}-\widetilde{C}(a)\frac{\omega^{-1}W_{0,1}(G(a),\omega^2\lambda)}{W_{0,2}(a,\lambda)}\nonumber\\
=&\frac{\omega[2i\omega^{-\frac{3}{2}}+o(1)]\exp\left[L(G^{2}(a),\omega^{2-m}\lambda)+L(G^{-1}(a),\omega^{-2}\lambda)\right]}{[2i\omega^{-\frac{2}{2}}+o(1)]\exp\left[L(G^{2}(a),\omega^{2-m}\lambda)+L(a,\lambda)\right]}\nonumber\\
&-\widetilde{C}(a)\frac{\omega^{-1}W_{0,1}(G(a),\omega^2\lambda)}{[2i\omega^{-1}+o(1)]\exp\left[L(G^{2}(a),\omega^{2-m}\lambda)+L(a,\lambda)\right]}\nonumber\\
=&[\omega^{\frac{1}{2}}+o(1)]\exp\left[L(G^{-1}(a),\omega^{-2}\lambda)-L(a,\lambda)\right]\nonumber\\
&-\frac{-\omega^{1+2\nu(a)}2\omega^{\mu(G(a))}}{[2i+o(1)]\exp\left[L(G^{2}(a),\omega^{2-m}\lambda)+L(a,\lambda)\right]}\nonumber
\end{align}
as $\lambda\to\infty$ such that
\begin{equation}\nonumber
-\pi+\delta\leq\pi-\frac{12\pi}{m+2}+\delta \leq \arg(\omega^{-2}\lambda)\leq \pi-\delta\quad\text{and}\quad
\pi-\frac{8\pi}{m+2}+\delta \leq \arg(\lambda)\leq \pi-\delta,
\end{equation}
that is,
\begin{equation}\nonumber
\pi-\frac{8\pi}{m+2}+\delta \leq \arg(\lambda)\leq \pi-\delta.
\end{equation}
Next, we use $2\nu(a)+\mu(G(a))=\mu(a)+4\nu(a)$ to get \eqref{asy_asy}.

Suppose $m=3$. Then $\omega^5=1$. Also, $W_{-3,0}(a,\lambda)=W_{2,0}(a,\lambda)$ since $f_{-3}(z,a,\lambda)=f_2(z,a,\lambda)$. Thus, we have that
\begin{align}
C(a,\lambda)=&\frac{W_{-1,2}(a,\lambda)}{W_{0,2}(a,\lambda)}-\widetilde{C}(a)\frac{W_{1,2}(a,\lambda)}{W_{0,2}(a,\lambda)}\nonumber\\
=&\frac{\omega^{-2} W_{-3,0}(G^{2}(a),\omega^{4}\lambda)}{W_{0,2}(a,\lambda)}-\widetilde{C}(a)\frac{\omega^{-1}W_{0,1}(G(a),\omega^2\lambda)}{W_{0,2}(a,\lambda)}\nonumber\\
=&-\frac{\omega^{-2} W_{0,2}(G^{2}(a),\omega^{4}\lambda)}{W_{0,2}(a,\lambda)}-\widetilde{C}(a)\frac{\omega^{-1}W_{0,1}(G(a),\omega^2\lambda)}{W_{0,2}(a,\lambda)}\nonumber\\
=&-\frac{\omega^{-2} W_{0,2}(G^{2}(a),\omega^{-1}\lambda)}{W_{0,2}(a,\lambda)}-\widetilde{C}(a)\frac{\omega^{-1}W_{0,1}(G(a),\omega^2\lambda)}{W_{0,2}(a,\lambda)}\nonumber\\
=&-\frac{\omega^{-2}[2i\omega^{-\frac{2}{2}}+o(1)]\exp\left[L(G^{4}(a),\omega^{-2}\lambda)+L(G^{2}(a),\omega^{-1}\lambda)\right]}{[2i\omega^{-\frac{2}{2}}+o(1)]\exp\left[L(G^{2}(a),\omega^{-1}\lambda)+L(a,\lambda)\right]}\nonumber\\
&-\widetilde{C}(a)\frac{\omega^{-1}W_{0,1}(G(a),\omega^2\lambda)}{[2i\omega^{-1}+o(1)]\exp\left[L(G^{2}(a),\omega^{-1}\lambda)+L(a,\lambda)\right]}\nonumber\\
=&[-\omega^{-2}+o(1)]\exp\left[L(G^{4}(a),\omega^{-2}\lambda)-L(a,\lambda)\right]\nonumber\\
&-\frac{-\omega^{1+2\nu(a)}2\omega^{\mu(G(a))}}{[2i+o(1)]\exp\left[L(G^{2}(a),\omega^{2-m}\lambda)+L(a,\lambda)\right]},\nonumber
\end{align}
as $\lambda\to\infty$ such that
\begin{equation}\nonumber
-\pi+\delta\leq\pi-\frac{8\pi}{5}+\delta \leq \arg(\omega^{-1}\lambda)\leq \pi-\delta\quad\text{and}\quad
\pi-\frac{8\pi}{5}+\delta \leq \arg(\lambda)\leq \pi-\delta,
\end{equation}
that is,
\begin{equation}\nonumber
\pi-\frac{6\pi}{5}+\delta \leq \arg(\lambda)\leq \pi-\delta.
\end{equation}

In order to show the last assertion, we suppose that  $C(a,\lambda)=0$ for some $\lambda$ in \eqref{sector2} with large $|\lambda|$. Then from the asymptotic expansion \eqref{asy_asy}, we have  
\begin{equation}\label{blows_up}
\exp\left[L(G^{-1}(a),\omega^{-2}\lambda)+L(G^{2}(a),\omega^{2-m}\lambda)\right]=[i\omega^{-\frac{1}{2}+\mu(a)+4\nu(a)}+o(1)].
\end{equation}
 By Corollary \ref{lemma_decay},  
\begin{align}
\Re\left(L(G^{-1}(a),\omega^{-2}\lambda)\right)&=K_{m}\cos\left(\frac{m+2}{2m}\arg(\omega^{-2}\lambda)\right)|\lambda|^{\frac{1}{2}+\frac{1}{m}}(1+o(1))\nonumber\\
&=K_{m}\cos\left(-\frac{2\pi}{m}+\frac{m+2}{2m}\arg(\lambda)\right)|\lambda|^{\frac{1}{2}+\frac{1}{m}}(1+o(1)),\nonumber\\
\Re\left(L(G^{2}(a),\omega^{2-m}\lambda)\right)&=K_{m}\cos\left(\frac{m+2}{2m}\arg(\omega^{2-m}\lambda)\right)|\lambda|^{\frac{1}{2}+\frac{1}{m}}(1+o(1))\nonumber\\
&=-K_{m}\cos\left(\frac{2\pi}{m}+\frac{m+2}{2m}\arg(\lambda)\right)|\lambda|^{\frac{1}{2}+\frac{1}{m}}(1+o(1)).\nonumber
\end{align}
Note that if $m\geq 6$, then $0<\delta\leq\arg(\lambda)\leq\pi-\delta$ in \eqref{sector2}.
Since
\begin{align}
&\cos\left(-\frac{2\pi}{m}+\frac{m+2}{2m}\arg(\lambda)\right)-\cos\left(\frac{2\pi}{m}+\frac{m+2}{2m}\arg(\lambda)\right)\nonumber\\
&=2\sin\left(\frac{2\pi}{m}\right)\sin\left(\frac{m+2}{2m}\arg(\lambda)\right)>0,\nonumber
\end{align}
we see that
$$\Re\left(L(G^{-1}(a),\omega^{-2}\lambda)+L(G^{2}(a),\omega^{2-m}\lambda)\right)\to+\infty,$$
as $\lambda\to\infty$ in \eqref{sector2}, and hence the left hand side of \eqref{blows_up} blows up. Thus, $C(a,\lambda)$
 cannot have infinitely many zeros in \eqref{sector2}. This completes the proof.
\end{proof}

The next theorem covers a sector in the lower half plane, connecting sectors \eqref{sector1} and \eqref{sector4}.
\begin{theorem}
Suppose that $a\in\C^{m-1}$ and $0<\delta<\frac{\pi}{m+2}$. If $m\geq 4$ then
\begin{align}
C(a,\lambda)=&\left[-i\omega^{1+\mu(a)}+o(1)\right]\exp\left[-L(a,\omega^{-m-2}\lambda)-L(G^{-2}(a),\omega^{-4}\lambda)\right]\nonumber\\
&+[\omega^{\frac{1}{2}+2\nu(a)}+o(1)]\exp\left[L(G(a),\omega^{-m}\lambda)-L(a,\omega^{-m-2}\lambda)\right],\nonumber
\end{align}
as $\lambda\to\infty$ in the sector
\begin{equation}\label{sector3}
\pi+\delta \leq \arg(\lambda)\leq \pi+\frac{8\pi}{m+2}-\delta.
\end{equation}
If $m=3$ then 
\begin{align}
C(a,\lambda)=&[-i\omega^{\frac{7}{4}}+o(1)]\exp\left[-L(a,\omega^{-5}\lambda)-L(G^{-2}(a),\omega^{-4}\lambda)\right]\nonumber\\
&+[\omega^3+o(1)]\exp\left[L(G(a),\omega^{-3}\lambda)-L(a,\omega^{-5}\lambda)\right],\nonumber
\end{align}
as $\lambda\to\infty$ in the sector 
\begin{equation}\label{sector3-1}
\pi+\delta \leq \arg(\lambda)\leq 2\pi+\frac{\pi}{5}-\delta.
\end{equation}

Moreover, if $m\geq 6$ then there exists a constant $M_3>0$ such that $C(a,\lambda)\not=0$ for all $\lambda$  in the sector \eqref{sector3} if $|\lambda|\geq M_3$. 
\end{theorem}

\begin{proof}
Suppose that $m\geq 4$. Then from  \eqref{kplus1}, \eqref{kdef} and Lemma \ref{lemma7},
\begin{align}
C(a,\lambda)=&\frac{W_{-1,-2}(a,\lambda)}{W_{0,-2}(a,\lambda)}-\widetilde{C}(a)\frac{W_{1,-2}(a,\lambda)}{W_{0,-2}(a,\lambda)}\nonumber\\
=&\frac{W_{0,1}(G^{-2}(a),\omega^{-4}\lambda)}{W_{0,2}(G^{-2}(a),\omega^{-4}\lambda)}+\widetilde{C}(a)\frac{\omega^{-1}W_{0,-3}(G(a),\omega^{2}\lambda)}{\omega^2W_{0,2}(G^{-2}(a),\omega^{-4}\lambda)}\nonumber\\
=&\frac{W_{0,1}(G^{-2}(a),\omega^{-4}\lambda)}{[2i\omega^{-\frac{2}{2}}+o(1)]\exp\left[L(a,\omega^{-m-2}\lambda)+L(G^{-2}(a),\omega^{-4}\lambda)\right]}\nonumber\\
&-\widetilde{C}(a)\frac{[2i\omega^{-\frac{3}{2}}+o(1)]\exp\left[L(G(a),\omega^{-m}\lambda)+L(G^{-2}(a),\omega^{-4}\lambda)\right]}{[2i\omega^{-\frac{2}{2}}+o(1)]\exp\left[L(a,\omega^{-m-2}\lambda)+L(G^{-2}(a),\omega^{-4}\lambda)\right]}\nonumber\\
=&\frac{2\omega^{\mu(G^{-2}(a))}}{[2i\omega^{-1}+o(1)]\exp\left[L(a,\omega^{-m-2}\lambda)+L(G^{-2}(a),\omega^{-4}\lambda)\right]}\nonumber\\
&+[\omega^{-\frac{1}{2}}+o(1)]\omega^{1+2\nu(a)}\frac{\exp\left[L(G(a),\omega^{-m}\lambda)+L(G^{-2}(a),\omega^{-4}\lambda)\right]}{\exp\left[L(a,\omega^{-m-2}\lambda)+L(G^{-2}(a),\omega^{-4}\lambda)\right]},\nonumber
\end{align}
as $\lambda\to \infty$ such that
\begin{equation}
\pi-\frac{12\pi}{m+2}+\delta\leq \arg(\omega^{-4}\lambda)\leq \pi-\delta\quad\text{and}\quad \pi-\frac{8\pi}{m+2}+\delta\leq \arg(\omega^{-4}\lambda)\leq \pi-\delta,
\end{equation}
that is,
\begin{equation}
\pi+\delta\leq\arg(\lambda)\leq\pi+\frac{8\pi}{m+2}-\delta,
\end{equation}
which is \eqref{sector3}.

Suppose that $m=3$. Then, 
\begin{align}
C(a,\lambda)=&\frac{W_{-1,-2}(a,\lambda)}{W_{0,-2}(a,\lambda)}-\widetilde{C}(a)\frac{W_{1,-2}(a,\lambda)}{W_{0,-2}(a,\lambda)}\nonumber\\
=&\frac{W_{0,1}(G^{-2}(a),\omega^{-4}\lambda)}{W_{0,2}(G^{-2}(a),\omega^{-4}\lambda)}+\widetilde{C}(a)\frac{\omega^{-1}W_{0,-3}(G(a),\omega^{2}\lambda)}{\omega^2W_{0,2}(G^{-2}(a),\omega^{-4}\lambda)}\nonumber\\
=&\frac{W_{0,1}(G^{-2}(a),\omega^{-4}\lambda)}{W_{0,2}(G^{-2}(a),\omega^{-4}\lambda)}-\omega^2\widetilde{C}(a)\frac{W_{0,2}(G(a),\omega^{-3}\lambda)}{W_{0,2}(G^{-2}(a),\omega^{-4}\lambda)}\nonumber\\
=&\frac{W_{0,1}(G^{-2}(a),\omega^{-4}\lambda)}{[2i\omega^{-\frac{2}{2}}+o(1)]\exp\left[L(a,\omega^{-5}\lambda)+L(G^{-2}(a),\omega^{-4}\lambda)\right]}\nonumber\\
&-\omega^2\widetilde{C}(a)\frac{[2i\omega^{-\frac{2}{2}}+o(1)]\exp\left[L(G^3(a),\omega^{-4}\lambda)+L(G(a),\omega^{-3}\lambda)\right]}{[2i\omega^{-\frac{2}{2}}+o(1)]\exp\left[L(a,\omega^{-5}\lambda)+L(G^{-2}(a),\omega^{-4}\lambda)\right]}\nonumber\\
=&\frac{2\omega^{\mu(G^{-2}(a))}}{[2i\omega^{-1}+o(1)]\exp\left[L(a,\omega^{-5}\lambda)+L(G^{-2}(a),\omega^{-4}\lambda)\right]}\nonumber\\
&+[\omega^{2}+o(1)]\omega^{1+2\nu(a)}\frac{\exp\left[L(G^{-2}(a),\omega^{-4}\lambda)+L(G(a),\omega^{-3}\lambda)\right]}{\exp\left[L(a,\omega^{-5}\lambda)+L(G^{-2}(a),\omega^{-4}\lambda)\right]},\nonumber
\end{align}
as $\lambda\to\infty$ such that
\begin{equation}\nonumber
\pi-\frac{8\pi}{5}+\delta \leq \arg(\omega^{-3}\lambda)\leq \pi-\delta\quad\text{and}\quad
\pi-\frac{8\pi}{5}+\delta \leq \arg(\omega^{-4}\lambda)\leq \pi-\delta,
\end{equation}
that is,
\begin{equation}\nonumber
\pi+\delta \leq \arg(\lambda)\leq \pi+\frac{6\pi}{5}-\delta.
\end{equation}
Finally,  the proof of the last assertion of this theorem follows as in the proof of Theorem  \ref{thm_sector2}.
\end{proof}

From the asymptotic expansions in the previous four theorems, one obtains the order of the entire function $\lambda\mapsto C(a,\lambda)$.
 The order of an entire function $g$ is defined by
$$\limsup_{r\rightarrow \infty}\frac{\log \log M(r,g)}{\log r},$$
where $M(r, g)=\max \{|g(re^{i\theta})|: 0\leq \theta\leq 2\pi\}$ for $r>0$.
If for some positive real numbers $\sigma,\, c_1,\, c_2$, we have $\exp[c_1 r^{\sigma}]\leq M(r,g)\leq  \exp[c_2 r^{\sigma}]$ for all large $r$, then the order of $g$ is $\sigma$. 
\begin{corollary}
The entire function $\lambda\mapsto C(a,\lambda)$ is of order $\frac{1}{2}+\frac{1}{m}$.
\end{corollary}
\begin{proof}
The sectors in \eqref{sector0}, \eqref{sector1}, \eqref{sector2}, \eqref{sector3}, cover a neighborhood of infinity in the complex plane. So the nonconstant entire function  $|C(a,\lambda)|$ is bounded above by $\exp\left[c_1|\lambda|^{\frac{1}{2}+\frac{1}{m}}\right]$ for some constant $c_1>0$. Also, along the ray $\arg(\lambda)=\pi$,  one can see from \eqref{re_part} and \eqref{asy_1} that $|C(a,\lambda)|$ is bounded below by $\exp\left[c_2|\lambda|^{\frac{1}{2}+\frac{1}{m}}\right]$ for some constant $c_2>0$. Hence, the order of $C(a,\cdot)$ is  $\frac{1}{2}+\frac{1}{m}$.
\end{proof}

\begin{remark}
{\rm Since the eigenvalues are the zeros of the entire function $\lambda\mapsto C(a,\lambda)$ of order $\frac{1}{2}+\frac{1}{m}\in (0,1)$, there are infinitely many discrete eigenvalues as was already mentioned in Theorem \ref{main2}. 
}
\end{remark}

\section{Asymptotic expansion of the eigenvalues: Proof of Theorem \ref{main_thm1}}\label{asymp_eigen}
In this section, we prove Theorem \ref{main_thm1} by using the asymptotic expansions of $C(a,\lambda)$ and $L(a,\lambda)$.

\begin{proof}[Proof of Theorem ~\ref{main_thm1}]
Recall that by Theorem \ref{iden_thm}, $\lambda$ is an eigenvalue of $H$ if and only if $C(a,\lambda)=0$.

For $m\geq 4$ and $a\in\C^{m-1}$ fixed, suppose that $C(a,\lambda)=0$ for some $\lambda$ with $|\lambda|$ large. Then from the asymptotic expansion of $C(a,\lambda)$ in Theorem \ref{zero_thm} we have
\begin{equation}\nonumber
[1+o(1)]\exp\left[L(G(a),\omega^2\lambda)-L(G^{-1}(a),\omega^{-2}\lambda)\right]
=-\omega^{-2\nu(a)},
\end{equation}
and absorbing $[1+o(1)]$ into the exponential function then  yields
\begin{equation}\nonumber
\exp\left[L(G(a),\omega^2\lambda)-L(G^{-1}(a),\omega^{-2}\lambda)+o(1)\right]
=-\omega^{-2\nu(a)}.
\end{equation}
Thus, from Lemma \ref{asy_lemma} if $m$ is odd, we infer
\begin{align}
&\ln\left(-\omega^{-2\nu(a)}\right)\nonumber\\
&=L(G(a),\omega^2\lambda)-L(G^{-1}(a),\omega^{-2}\lambda)+o(1)\nonumber\\
&=
\sum_{j=0}^{\lfloor\frac{m}{2}+1\rfloor}\left[ K_{m,j}(G(a))(\omega^{2}\lambda)^{\frac{1}{2}+\frac{1-j}{m}}-K_{m,j}(G^{-1}(a))(\omega^{-2}\lambda)^{\frac{1}{2}+\frac{1-j}{m}}\right]+o(1)\label{region_def}\\
&=2iK_{m,0}\sin\left(\frac{2\pi}{m}\right)\lambda^{\frac{1}{2}+\frac{1}{m}}+\sum_{j=1}^{\lfloor\frac{m}{2}+1\rfloor}c_{m,j}(a)\lambda^{\frac{1}{2}+\frac{1-j}{m}}+o(1),\nonumber
\end{align}
where
\begin{equation}\label{c_def}
c_{m,j}(a)= K_{m,j}(G(a))(\omega^{2})^{\frac{1}{2}+\frac{1-j}{m}}-K_{m,j}(G^{-1}(a))(\omega^{-2})^{\frac{1}{2}+\frac{1-j}{m}},\quad 1\leq j\leq \frac{m+1}{2}.
\end{equation}

Similarly, if $m$ is even, then from Lemma \ref{asy_lemma} we have \eqref{region_def} with $c_{m,j}(a)$ in \eqref{c_def} for $1\leq j\leq \frac{m}{2}$, and 
$$
c_{m,\frac{m}{2}+1}(a)= K_{m,\frac{m}{2}+1}(G(a))-K_{m,\frac{m}{2}+1}(G^{-1}(a))+\frac{b_{\frac{m}{2}+1}(a)}{m}\frac{8\pi i}{m+2},
$$
where we used $b_{\frac{m}{2}+1}(G^{-1}(a))=-b_{\frac{m}{2}+1}(a)=b_{\frac{m}{2}+1}(G(a)).$

Note that there exist  constants $M>0$ and $\varepsilon>0$ such that the function 
\begin{equation}\label{map_eq}
\lambda\mapsto L(G(a),\omega^2\lambda) -L(G^{-1}(a),\omega^{-2}\lambda)+o(1)
\end{equation} 
is continuous in the region $|\lambda|\geq M$ and $|\arg(\lambda)|\leq \varepsilon$.
From \eqref{region_def} we then see that the function \eqref{map_eq} maps the region $|\lambda|\geq M$ and $|\arg(\lambda)|\leq \varepsilon$ onto a region that contains the entire positive imaginary axis near infinity. 

Thus, from \eqref{region_def} we get that for every sufficiently large $n\in\N$  there exists $\lambda_n$ such that
\begin{equation}\nonumber
2iK_{m,0}\sin\left(\frac{2\pi}{m}\right)\lambda_n^{\frac{1}{2}+\frac{1}{m}}+\sum_{j=1}^{\lfloor\frac{m}{2}+1\rfloor}c_{m,j}(a)\lambda_n^{\frac{1}{2}+\frac{1-j}{m}}+o(1)=\left(2n+1-\frac{4\nu(a)}{m+2}\right)\pi i.
\end{equation}
Thus,
\begin{equation}\nonumber
\lambda_n^{\frac{1}{2}+\frac{1}{m}}+\sum_{j=1}^{\lfloor\frac{m}{2}+1\rfloor}\frac{c_{m,j}(a)}{2iK_{m,0}\sin\left(\frac{2\pi}{m}\right)}\lambda_n^{\frac{1}{2}+\frac{1-j}{m}}+o(1)=\frac{\left(2n+1-\frac{4\nu(a)}{m+2}\right)\pi}{2K_{m,0}\sin\left(\frac{2\pi}{m}\right)}.
\end{equation}
Let 
\begin{eqnarray}
d_{m,j}(a)=\left\{
              \begin{array}{rl}
              \frac{c_{m,j}(a)}{2iK_{m,0}\sin\left(\frac{2\pi}{m}\right)}  \quad & \text{if $1\leq j\leq \frac{m}{2}+1$,}\\
               \frac{c_{m,j}(a)+\frac{4\nu(a)}{m+2}\pi i}{2iK_{m,0}\sin\left(\frac{2\pi}{m}\right)} \quad & \text{if $m$ is even and $j=\frac{m}{2}+1$.}
              \end{array}
                         \right. \label{d_def}
\end{eqnarray}
Then
\begin{equation}\label{asy_eq3}
\lambda_n^{\frac{1}{2}+\frac{1}{m}}+\sum_{j=1}^{\lfloor\frac{m}{2}+1\rfloor}d_{m,j}(a)\lambda_n^{\frac{1}{2}+\frac{1-j}{m}}+o(1)=\frac{\left(2n+1\right)\pi}{2K_{m,0}\sin\left(\frac{2\pi}{m}\right)}.
\end{equation}
Introduce the decomposition
$
\lambda_n=\lambda_{0,n}+\lambda_{1,n},
$
where 
\begin{equation}\nonumber
 \lambda_{0,n}=\left(\frac{\left(2n+1\right)\pi}{2K_{m,0}\sin\left(\frac{2\pi}{m}\right)}\right)^{\frac{2m}{m+2}}\,\,\text{and}\quad \frac{\lambda_{1,n}}{\lambda_{0,n}}=o\left(1\right).
\end{equation}
Then from \eqref{asy_eq3} we have
\begin{align}
\lambda_{0,n}^{\frac{1}{2}+\frac{1}{m}}
&=\lambda_{0,n}^{\frac{1}{2}+\frac{1}{m}}\left(1+\frac{\lambda_{1,n}}{\lambda_{0,n}}\right)^{\frac{1}{2}+\frac{1}{m}}+\sum_{j=1}^{\lfloor\frac{m}{2}+1\rfloor}d_{m,j}(a)\lambda_{0,n}^{\frac{1}{2}+\frac{1-j}{m}}\left(1+\frac{\lambda_{1,n}}{\lambda_{0,n}}\right)^{\frac{1}{2}+\frac{1-j}{m}}+o(1)\nonumber\\
&=\lambda_{0,n}^{\frac{1}{2}+\frac{1}{m}}\left(1+\sum_{k=1}^{\infty}{\frac{1}{2}+\frac{1}{m}\choose k}\left(\frac{\lambda_{1,n}}{\lambda_{0,n}}\right)^k\right)\nonumber\\
&+\sum_{j=1}^{\lfloor\frac{m}{2}+1\rfloor}d_{m,j}(a)\lambda_{0,n}^{\frac{1}{2}+\frac{1-j}{m}}\left(1+\sum_{k=1}^{\infty}{\frac{1}{2}+\frac{1-j}{m}\choose k}\left(\frac{\lambda_{1,n}}{\lambda_{0,n}}\right)^k\right)+o(1).\nonumber
\end{align}
Thus,
\begin{align}
0&=\left(\frac{1}{2}+\frac{1}{m}\right)\frac{\lambda_{1,n}}{\lambda_{0,n}}+\sum_{k=2}^{\infty}{\frac{1}{2}+\frac{1}{m}\choose k}\left(\frac{\lambda_{1,n}}{\lambda_{0,n}}\right)^k\nonumber\\
&+\sum_{j=1}^{\lfloor\frac{m}{2}+1\rfloor}d_{m,j}(a)\lambda_{0,n}^{-\frac{j}{m}}\left(1+\sum_{k=1}^{\infty}{\frac{1}{2}+\frac{1-j}{m}\choose k}\left(\frac{\lambda_{1,n}}{\lambda_{0,n}}\right)^k\right)+o\left(\lambda_{0,n}^{-\frac{1}{2}-\frac{1}{m}}\right),\nonumber
\end{align}
and hence
\begin{align}
&\frac{\lambda_{1,n}}{\lambda_{0,n}}+\frac{2m}{m+2}\sum_{k=2}^{\infty}{\frac{1}{2}+\frac{1}{m}\choose k}\left(\frac{\lambda_{1,n}}{\lambda_{0,n}}\right)^k\nonumber\\
&+\frac{2m}{m+2}\sum_{j=1}^{\lfloor\frac{m}{2}+1\rfloor}d_{m,j}(a)\lambda_{0,n}^{-\frac{j}{m}}\left(\sum_{k=1}^{\infty}{\frac{1}{2}+\frac{1-j}{m}\choose k}\left(\frac{\lambda_{1,n}}{\lambda_{0,n}}\right)^k\right)+o\left(\lambda_{0,n}^{-\frac{1}{2}-\frac{1}{m}}\right)\nonumber\\
&=-\frac{2m}{m+2}\sum_{j=1}^{\lfloor\frac{m}{2}+1\rfloor}d_{m,j}(a)\lambda_{0,n}^{-\frac{j}{m}}.\label{asy_eq4}
\end{align}
Thus, one concludes
$
\frac{\lambda_{1,n}}{\lambda_{0,n}}=\lambda_{2,n}+\lambda_{3,n},
$
where
\begin{equation}\label{ex_eq1}
\lambda_{2,n}=-\frac{2m}{m+2}d_{m,1}(a)\lambda_{0,n}^{-\frac{1}{m}}\,\,\text{ and }\,\,\lambda_{3,n}=o\left(\lambda_{0,n}^{-\frac{1}{m}}\right).
\end{equation}
Next, from  \eqref{ex_eq1} along with \eqref{asy_eq4} we have
\begin{align}
&\lambda_{2,n}+\lambda_{3,n}+\frac{2m}{m+2}\sum_{k=2}^{\infty}{\frac{1}{2}+\frac{1}{m}\choose k}\left(\lambda_{2,n}+\lambda_{3,n}\right)^k\nonumber\\
&+\frac{2m}{m+2}\sum_{j=1}^{\lfloor\frac{m}{2}+1\rfloor}d_{m,j}(a)\lambda_{0,n}^{-\frac{j}{m}}\left(\sum_{k=1}^{\infty}{\frac{1}{2}+\frac{1-j}{m}\choose k}\left(\lambda_{2,n}+\lambda_{3,n}\right)^k\right)+o\left(\lambda_{0,n}^{-\frac{1}{2}-\frac{1}{m}}\right)\nonumber\\
&=-\frac{2m}{m+2}\sum_{j=1}^{\lfloor\frac{m}{2}+1\rfloor}d_{m,j}(a)\lambda_{0,n}^{-\frac{j}{m}}.\label{asy_eq5}
\end{align}
Thus,
\begin{align}
&\lambda_{3,n}+\frac{2m}{m+2}\sum_{k=2}^{\infty}{\frac{1}{2}+\frac{1}{m}\choose k}\sum_{\ell=0}^k{k\choose \ell}\lambda_{2,n}^{\ell}\lambda_{3,n}^{k-\ell}\nonumber\\
&+\frac{2m}{m+2}\sum_{j=1}^{\lfloor\frac{m}{2}+1\rfloor}d_{m,j}(a)\lambda_{0,n}^{-\frac{j}{m}}\left(\sum_{k=1}^{\infty}{\frac{1}{2}+\frac{1-j}{m}\choose k}\sum_{\ell=0}^k{k\choose \ell}\lambda_{2,n}^{\ell}\lambda_{3,n}^{k-\ell}\right)+o\left(\lambda_{0,n}^{-\frac{1}{2}-\frac{1}{m}}\right)\nonumber\\
&=-\frac{2m}{m+2}\sum_{j=2}^{\lfloor\frac{m}{2}+1\rfloor}d_{m,j}(a)\lambda_{0,n}^{-\frac{j}{m}},\label{asy_eq6}
\end{align}
and hence
\begin{align}
&\lambda_{3,n}+\frac{2m}{m+2}\sum_{k=2}^{\infty}{\frac{1}{2}+\frac{1}{m}\choose k}\sum_{\ell=0}^{k-1}{k\choose \ell}\lambda_{2,n}^{\ell}\lambda_{3,n}^{k-\ell}\nonumber\\
&+\frac{2m}{m+2}\sum_{j=1}^{\lfloor\frac{m}{2}+1\rfloor}d_{m,j}(a)\lambda_{0,n}^{-\frac{j}{m}}\left(\sum_{k=1}^{\infty}{\frac{1}{2}+\frac{1-j}{m}\choose k}\sum_{\ell=0}^{k-1}{k\choose \ell}\lambda_{2,n}^{\ell}\lambda_{3,n}^{k-\ell}\right)+o\left(\lambda_{0,n}^{-\frac{1}{2}-\frac{1}{m}}\right)\nonumber\\
&=-\frac{2m}{m+2}\sum_{j=2}^{\lfloor\frac{m}{2}+1\rfloor}d_{m,j}(a)\lambda_{0,n}^{-\frac{j}{m}}-\frac{2m}{m+2}\sum_{k=2}^{\infty}{\frac{1}{2}+\frac{1}{m}\choose k}\lambda_{2,n}^{k}\nonumber\\
&-\frac{2m}{m+2}\sum_{j=1}^{\lfloor\frac{m}{2}+1\rfloor}d_{m,j}(a)\lambda_{0,n}^{-\frac{j}{m}}\left(\sum_{k=1}^{\infty}{\frac{1}{2}+\frac{1-j}{m}\choose k}\lambda_{2,n}^{k}\right).\label{asy_eq7}
\end{align}
So we choose 
\begin{equation}\label{eq_asy}
\lambda_{3,n}=\lambda_{4,n}+\lambda_{5,n},
\end{equation}
where
\begin{align}
\lambda_{4,n}&=-\frac{2m}{m+2}d_{m,2}(a)\lambda_{0,n}^{-\frac{2}{m}}-\frac{2m}{m+2}{\frac{1}{2}+\frac{1}{m}\choose 2}\lambda_{2,n}^{2}-\frac{m\,d_{m,1}(a)}{m+2}\lambda_{0,n}^{-\frac{1}{m}}\lambda_{2,n}\nonumber\\
&=\left(-\frac{2m}{m+2}d_{m,2}(a)+\left(\frac{2m^2}{(m+2)^2}-\left(\frac{2m}{m+2}\right)^3{\frac{1}{2}+\frac{1}{m}\choose 2}\right)d_{m,1}(a)^2\right)\lambda_{0,n}^{-\frac{2}{m}},\nonumber\\
\lambda_{5,n}&=o\left(\lambda_{0,n}^{-\frac{2}{m}}\right).\nonumber
\end{align}
Next, we replace $\lambda_{3,n}$ in \eqref{asy_eq7} by \eqref{eq_asy}. Upon iterating this process we get
\begin{align}
\lambda_n&=\lambda_{0,n}+\lambda_{1,n}=\lambda_{0,n}\left(1+\frac{\lambda_{1,n}}{\lambda_{0,n}}\right)\nonumber\\
&=\lambda_{0,n}\left(1+\lambda_{2,n}+\lambda_{3,n}\right)\nonumber\\
&=\lambda_{0,n}\left(1+\lambda_{2,n}+\lambda_{4,n}+\lambda_{5,n}\right)\nonumber\\
&\cdots\nonumber\\
&=\lambda_{0,n}\left(1+\sum_{\ell=1}^{\lfloor\frac{m}{2}+1\rfloor}e_{\ell}(a)\lambda_{0,n}^{-\frac{\ell}{m}}+o\left(\lambda_{0,n}^{-\frac{1}{2}-\frac{1}{m}}\right)\right),\label{e_def}
\end{align}  
as $n\to+\infty$, that is, \eqref{asym_eq1}.

Suppose that $m=3$. For this case we will use the asymptotic expansion in Theorem \ref{thm_sector2} that is valid in \eqref{sector2-1}. Similarly to what we did for the case $m\geq 4$, if $C(a,\lambda)=0$ then from the asymptotic expansion in Theorem \ref{thm_sector2} 
we have
\begin{equation}\nonumber
[1+o(1)]\exp\left[L(G^{4}(a),\omega^{-2}\lambda)+L(G^{2}(a),\omega^{-1}\lambda)\right]=-i\omega^{\frac{15}{4}}.
\end{equation} Thus, 
since 
$L(a,\lambda)=K_{3,0}(a)\lambda^{\frac{5}{6}}+K_{3,1}(a)\lambda^{\frac{3}{6}}+K_{3,2}(a)\lambda^{\frac{1}{6}}+o(1),$ we have
\begin{align}
&L(G^{4}(a),\omega^{-2}\lambda)+L(G^{2}(a),\omega^{-1}\lambda)+o(1)\nonumber\\
&=K_{3,0}(G^{4}(a))\left(\omega^{-2}\lambda\right)^{\frac{5}{6}}+K_{3,1}(G^{4}(a))\left(\omega^{-2}\lambda\right)^{\frac{3}{6}}+K_{3,2}(G^{4}(a))\left(\omega^{-2}\lambda\right)^{\frac{1}{6}}\nonumber\\
&+K_{3,0}(G^{2}(a))\left(\omega^{-1}\lambda\right)^{\frac{5}{6}}+K_{3,1}(G^{2}(a))\left(\omega^{-1}\lambda\right)^{\frac{3}{6}}+K_{3,2}(G^{2}(a))\left(\omega^{-1}\lambda\right)^{\frac{1}{6}}+o(1)\nonumber\\
&=K_{3,0}\left(e^{-i\frac{2\pi}{3}}+e^{-i\frac{\pi}{3}}\right)\lambda^{\frac{5}{6}}+c_{3,1}(a)\lambda^{\frac{3}{6}}+c_{3,2}(a)\lambda^{\frac{1}{6}}+o(1)\nonumber\\
&=-2iK_{3,0}\sin\left(\frac{2\pi}{3}\right)\lambda^{\frac{5}{6}}+c_{3,1}(a)\lambda^{\frac{3}{6}}+c_{3,2}(a)\lambda^{\frac{1}{6}}+o(1).\nonumber
\end{align}
So the continuous function $\lambda\mapsto L(G^{4}(a),\omega^{-2}\lambda)+L(G^{2}(a),\omega^{-1}\lambda)+o(1)$ maps a neighborhood of the positive real axis near infinity onto a neighborhood of the negative imaginary  axis near infinity. Hence, there exist a sequence of $\lambda_n$ near the positive real axis such that for all large enough positive integers $n$,
\begin{equation}\nonumber
-2iK_{3,0}\sin\left(\frac{2\pi}{3}\right)\lambda_n^{\frac{5}{6}}+c_{3,1}(a)\lambda_n^{\frac{3}{6}}+c_{3,2}(a)\lambda_n^{\frac{1}{6}}+o(1)=\ln\left(-i\omega^{\frac{15}{4}}\right)=\left(\pi-2(n+1)\pi\right)i.
\end{equation}   
From this result one concludes that the asymptotic expansion \eqref{asym_eq1} holds for $m=3$ as well similarly to the proof for the case $m\geq 4$.
\end{proof}

\section{Proof of Theorem \ref{ineq_eq}}
\begin{proof}[Proof of Theorem ~\ref{ineq_eq}]
First, note from \eqref{asym_eq1} that $\arg(\lambda_n)\to 0$ as $n\to+\infty$.

Next, we have
\begin{align}
\lambda_{0,n+1}&=\left(\frac{\left(2n+3\right)\pi}{2K_{m,0}\sin\left(\frac{2\pi}{m}\right)}\right)^{\frac{2m}{m+2}}\nonumber\\
&=\left(\frac{\left(2n+1\right)\pi}{2K_{m,0}\sin\left(\frac{2\pi}{m}\right)}+\frac{2\pi}{2K_{m,0}\sin\left(\frac{2\pi}{m}\right)}\right)^{\frac{2m}{m+2}}\nonumber\\
&=\lambda_{0,n}\left(1+\frac{2}{2n+1}\right)^{\frac{2m}{m+2}}\nonumber\\
&=\lambda_{0,n}\left(1+\frac{2m}{m+2}\frac{2}{2n+1}+O\left(\frac{1}{n^2}\right)\right)\nonumber\\
&=\lambda_{0,n}+\frac{2m\pi}{(m+2)K_{m,0}\sin\left(\frac{2\pi}{m}\right)}\lambda_{0,n}^{1-\frac{1}{2}-\frac{1}{m}}+o\left(\lambda_{0,n}^{\frac{1}{2}-\frac{1}{m}}\right).\label{ex_eq2}
\end{align}
Thus,
\begin{equation}\nonumber
\lambda_{n+1}-\lambda_{n}\underset{n\to+\infty}{=}\frac{2m\pi}{(m+2)K_{m,0}\sin\left(\frac{2\pi}{m}\right)}\lambda_{0,n}^{\frac{1}{2}-\frac{1}{m}}+o\left(\lambda_{0,n}^{\frac{1}{2}-\frac{1}{m}}\right),
\end{equation}
and hence, $|\lambda_{n+1}-\lambda_n|\to\infty$ and $\arg(\lambda_{n+1}-\lambda_n)\to 0$  as $n\to+\infty$.
Since $\arg(\lambda_n)\to 0$ (and $\arg(\lambda_{n+1})\to 0$) as $n\to+\infty$,  there exists $N\in\N$ such that $|\lambda_n|<|\lambda_{n+1}|$ if $n\geq N$.
\end{proof}

\begin{remark}
{\em Here we will show that if $a\in\R^{m-1}$, then $e_{\ell}(a)\in\R$ for all $1\leq \ell\leq \frac{m}{2}+1$ with $e_{\ell}(a)$ defined in \eqref{e_def}.

From \eqref{def_K} one can see that $\overline{K_{m,j}(G^{-1}(\overline{a}))}=K_{m,j}(G(a))$. Next, suppose that $a\in\R^{m-1}$.
If $m\geq 4$ then from \eqref{c_def}, 
\begin{align}
ic_{m,j}(a)&=i\left( K_{m,j}(G(a))(\omega^{2})^{\frac{1}{2}+\frac{1-j}{m}}-K_{m,j}(G^{-1}(a))(\omega^{-2})^{\frac{1}{2}+\frac{1-j}{m}}\right)\nonumber\\
&=i\left( K_{m,j}(G(a))(\omega^{2})^{\frac{1}{2}+\frac{1-j}{m}}- \overline{K_{m,j}(G(a))(\omega^{2})^{\frac{1}{2}+\frac{1-j}{m}}}\right)\in\R,\quad 1\leq j\leq \frac{m}{2}+1.\nonumber
\end{align}
So by \eqref{d_def}, $d_{m,j}(a)\in\R$ for all $1\leq j\leq \frac{m}{2}+1,$ and hence by \eqref{e_def}, $e_{\ell}(a)\in\R$ for all $1\leq \ell\leq \frac{m}{2}+1$. 

If $m=3$ then one can show $e_{\ell}(a)\in\R$ for $\ell=1,\,2$, using the formulas at the end of the proof of Theorem \ref{main_thm1}.
}
\end{remark}

\subsection*{{\bf Acknowledgments}}
 The author thanks Fritz Gesztesy and Richard Laugesen for critical reading of this manuscript and suggestions.

{\sc email contact:}  kcshin@math.missouri.edu

\begin{thebibliography}{10}

\bibitem{Bender-1}
C. M. Bender, M. Berry, P. N. Meisinger, V. M. Savage and M. Simsek
\newblock Complex WKB Analysis of Energy-Level Degeneracies of Non-Hermitian Hamiltonians,
\newblock{\em J. Phys. A: Math. Gen.}, 34:L31--L36, 2001.

\bibitem{Bender}
C. M. Bender and S. Boettcher.
\newblock Real spectra in non-Hermitian Hamiltonians having $\mathcal{PT}$-symmetry.
\newblock {\em Phys. Rev. Lett.}, 80:5243--5246, 1998. 

\bibitem{Bender2}
C. M. Bender and A. Turbiner.
\newblock Analytic continuation of eigenvalue problems.
\newblock {\em Phys. Lett. A}, 173:442--446, 1993.

\bibitem{Bender5}
C. M. Bender and E. J. Weniger. 
\newblock Numerical evidence that the perturbation expansion for a non-Hermitian $\mathcal{PT}$-symmetric Hamiltonian is Stieltjes.
\newblock {\em J. Math. Phys.}, 42:2167--2183, 2001.

\bibitem{CGM}
E. Caliceti, S. Graffi and M. Maioli.
\newblock Perturbation theory of odd anharmonic oscillators.
\newblock {\em Comm. Math. Phys.}, 75:51--66, 1980.


\bibitem{Pham}
E. Delabaere and F. Pham.
\newblock Eigenvalues of complex Hamiltonians with $\mathcal{PT}$-symmetry I, II.
\newblock {\em Phys. Lett. A}, 250:25--32, 1998.

\bibitem{Delabaere}
E. Delabaere and D. T. Trinh.
\newblock Spectral analysis of the complex cubic oscillator.
\newblock {\em J. Phys. A: Math. Gen.}, 33:8771--8796, 2000. 


\bibitem{Dorey}
P. Dorey, C. Dunning and R. Tateo.
\newblock Spectral equivalences, Bethe ansatz equations, and reality properties in $\mathcal{PT}$-symmetric quantum mechanics.
\newblock {\em J. Phys. A: Math. Gen}, 34:5679--5704, 2001.


\bibitem{Dorey2}
P. Dorey and R. Tateo.
\newblock On the relation between Stokes multipliers and $T$-$Q$ systems of conformal field theory.
\newblock {\em Nucl. Phys. B}, 563:573--602, 1999.


\bibitem{Handy2}
C. R. Handy. 
\newblock Generating converging bounds to the (complex) discrete states of the $P^2 + iX^3 + i\alpha X$ Hamiltonian.
\newblock {\em J. Phys. A: Math. Gen.}, 34:5065--5081, 2001.

\bibitem{Handy1}
C. R. Handy, D. Khan, Xiao-Qian Wang and C. J. Tymczak.
\newblock Multiscale reference function analysis of the  $\mathcal{PT}$ symmetry breaking solutions for the $P^2 + iX^3 + i\alpha X$ Hamiltonian.
\newblock {\em J. Phys. A: Math. Gen.}, 34:5593--5602, 2001.


\bibitem{Hille}
E. Hille.
\newblock {\em Lectures on Ordinary Differential Equations}.
\newblock Addison-Wesley, Reading, Massachusetts, 1969.


\bibitem{mez}
G. A. Mezincescu.
\newblock Some properties of eigenvalues and eigenfunctions of the cubic oscillator with imaginary coupling constant.
\newblock {\em J. Phys. A: Math. Gen.}, 33:4911--4916, 2000.

\bibitem{Ali1}
A. Mostafazadeh.
\newblock Pseudo-Hermiticity versus PT Symmetry: The necessary condition for the reality of the spectrum of a non-Hermitian Hamiltonian.
\newblock {\em J. Math. Phys.}, 43:205--214, 2002.


\bibitem{Shin1}
K. C. Shin.
\newblock On the eigenproblems of $\mathcal{PT}$-symmetric oscillators.
\newblock {\em J. Math. Phys.}, 42:2513--2530, 2001.

\bibitem{Shin}
K. C. Shin.
\newblock On the reality of the eigenvalues for a class of $\mathcal{PT}$-symmetric oscillators,
\newblock {\em Comm. Math. Phys.}, 229(3):543--564, 2002.

\bibitem{Simon}
B. Simon.
\newblock Coupling Constant Analyticity for the Anharmonic Oscillator.
\newblock {\em Ann. Phys.}, 58:76--136, 1970.




\bibitem{Sibuya}
Y. Sibuya.
\newblock {\em Global theory of a second order linear ordinary differential equation with a polynomial coefficient}.
\newblock North-Holland Publishing Company, Amsterdam-Oxford, 1975.


\bibitem{Znojil}
M. Znojil.
\newblock Spiked and  $\mathcal{PT}$-symmetrized decadic potentials supporting elementary $N$-plets of bound states.
\newblock {\em J. Phys. A: Math. Gen.}, 33:4911--4916, 2000. 

\end{thebibliography}
\end{document}